\documentclass[conference, 10pt]{IEEEtran}
\usepackage[pdftex]{graphicx}
\usepackage[cmex10]{amsmath}
\usepackage{amsthm}
\usepackage{amssymb}
\usepackage[all,2cell]{xy}\UseAllTwocells\SilentMatrices
\theoremstyle{plain}

\newtheorem{thm}{Theorem}

\newtheorem{lem}[thm]{Lemma}

\theoremstyle{definition}
\newtheorem{df}[thm]{Definition}
\newtheorem{rem}[thm]{Remark}
\newtheorem{eg}[thm]{Example}

\DeclareMathOperator*{\id}{id} 
\DeclareMathOperator*{\pr}{pr} 
\def\attach{\rightsquigarrow} 
\def\shf{\mathcal}            

\begin{document}
\title{Understanding networks and their behaviors using sheaf theory}

\author{\IEEEauthorblockN{Michael Robinson}
\IEEEauthorblockA{Mathematics and Statistics\\
American University\\
Washington, DC, USA\\
michaelr@american.edu}\\}

\IEEEspecialpapernotice{(Invited Paper)}

\maketitle

\begin{abstract}
Many complicated network problems can be easily understood on small networks.  Difficulties arise when small networks are combined into larger ones. Fortunately, the mathematical theory of sheaves was constructed to address just this kind of situation; it extends locally-defined structures to globally valid inferences by way of consistency relations. This paper exhibits examples in network monitoring and filter hardware where sheaves have useful descriptive power. 
\end{abstract}


\section{Introduction}

This paper serves as an invitation to the signal processing community to think about how consistency between pieces of local information serves to drive global inferences.  Because \emph{sheaves} are the mathematical tool for manipulating local data, they could be useful for solving certain signal processing problems.  For example, this article shows that a sampling theory for sheaves can lead to an inference methodology for monitoring water distribution networks.  In a rather different context, we show that sheaves are a natural model of signal processing hardware.
 
Sheaves were introduced by Leray during World War II to study partial differential equations.  They remained enshrouded within ``pure'' mathematics, where they enjoyed a central role in Grothendieck's algebraic geometry programme.  (See the foreword by C. Houzel in \cite{Kashiwara_1990} for a clear discussion of the early history of sheaves.)  An early departure from this abstract setting to a more concrete, computational setting occurred in \cite{Baclawski_1975}, which treated sheaves over posets.  However, cellular sheaves introduced by Shepard \cite{Shepard_1980}, are usually more natural for engineering contexts.  These ideas lay dormant until a recent flurry of activity \cite{Lilius_1993,ghrist_2011,RobinsonLogic,robinson_2013}.

\section{The main ideas of sheaf theory} 
A sheaf is a mathematical datastructure for storing local information over a topological space.  This article addresses sheaves over \emph{abstract simplicial complexes}, for which the theory is minimally complicated.  (The interested reader is encouraged to explore \cite{Curry} for a more extensive -- but readable -- introduction to sheaves over cellular spaces.)

\begin{df}
An \emph{abstract simplicial complex} over a set $A$ is a collection $X$ of (possibly ordered) subsets of $A$, for which $x \in X$ implies that every subset of $x$ is also in $X$.  We call each $x\in X$ with $k+1$ elements a \emph{$k$-face} of $X$, referring to the number $k$ as its \emph{dimension}.  Zero dimensional faces (singleton subsets of $A$) are called \emph{vertices}, and one dimensional faces are called \emph{edges}.  We say that a face $a$ \emph{includes into} a face $b$ (written $a\attach b$) whenever $a$ is a proper subset of $b$.  
\end{df}

\begin{eg}
A graph $G=(V,E)$ with vertices $V=\{v_1,\dotsc\}$ and edges $E=\{e_1,\dotsc\}$ can be given the structure of a simplicial complex $X=\{\{v_1\},\dotsc,e_1,\dotsc\}$.  Observe that if $\{v_1,v_2\}$ is an edge in $X$, then $\{v_1\},\{v_2\}\in X$ automatically by this construction.
\end{eg}

Given a simplicial complex, a sheaf is merely the assignment of vector-valued data to each face that is compatible with the inclusions of faces.

\begin{df}
A \emph{sheaf} $\shf{F}$ on an abstract simplicial complex $X$ consists of the assignment of 
\begin{enumerate}
\item a vector space $\shf{F}(a)$ to each face $a$ of $X$ (called the \emph{stalk} at $a$), and 
\item a linear map $\shf{F}(a\attach b):\shf{F}(a)\to \shf{F}(b)$ (called the \emph{restriction} along $a\attach b$) to each inclusion of faces $a\attach b$, so that
\item $\shf{F}(b\attach c) \circ \shf{F}(a\attach b)=\shf{F}(a\attach c)$ whenever $a\attach b \attach c$.
\end{enumerate}
\end{df}

\begin{df}
Suppose $X$ is the abstract simplicial complex model of $\mathbb{R}$ whose vertices are given by the set of integers and whose edges are given by pairs $\{n,n+1\}$.  The \emph{$N$-term grouping sheaf $\shf{V}^{(N)}$} is given by the diagram written over $X$
\begin{equation*}
\xymatrix{
\ar[r]&V^{N-1}&V^N\ar[r]^{\sigma_+}\ar[l]^{\sigma_-}&V^{N-1}&V^N\ar[r]^{\sigma_+}\ar[l]^{\sigma_-}&
}
\end{equation*}
where $\sigma_\pm$ are the $(N-1)\times N$ matrices
\begin{equation*}
\sigma_-=\begin{pmatrix}
0 & I_{(N-1)\times(N-1)}
\end{pmatrix}
\text{ and }
\sigma_+=\begin{pmatrix}
I_{(N-1)\times(N-1)} & 0
\end{pmatrix}.
\end{equation*}
\end{df}

The $N$-term grouping sheaf represents the time evolution of the contents of a $N$-word shift register.  The underlying simplicial complex structure for this sheaf represents a discrete timeline, with each vertex representing a timestep.  The stalk of $\shf{V}^{(N)}$ at each vertex is $V^N$, which represents the contents of the shift register at that timestep.

\begin{eg}
To see how this works, consider a 3-word shift register that stores integers.  Suppose that at time 0 it stores the vector (1,1,9), and at time 1 it loads a 2 into the last slot so that it contains (1,9,2).  Observe that the $\sigma_\pm$ maps defined above compute what is preserved during the transition from time 0 to time 1: $\sigma_+ (1,9,2) = (1,9) = \sigma_- (1,1,9)$.  
\end{eg}

The example shows that local consistency between values from the stalks of a sheaf leads to information that is consistent more globally.  This motivates the definition of a \emph{section} of a sheaf.

\begin{df}
Suppose that $\shf{F}$ is a sheaf on an abstract simplicial complex $X$.  A \emph{global section} $s$ assigns a value $s(a)\in \shf{F}(a)$ to each $a\in X$ so that for each inclusion $a\attach b$ of faces, $\shf{F}(a\attach b)s(a)=s(b)$.  The set of global sections forms a vector space $\Gamma \shf{F}$.
\end{df}

Although sheaves are useful descriptors, it is more important to focus on \emph{sheaf morphisms}, which are ways to manipulate sheaves.  
\begin{df}
A \emph{sheaf morphism} $f:\shf{F}\to \shf{G}$ of sheaves on an abstract simplicial complex $X$ assigns a linear map $f_a:\shf{F}(a) \to \shf{G}(a)$ to each face $a$ so that for every inclusion $a\attach b$ of faces of $X$, $f_b \circ \shf{F}(a\attach b) = \shf{G}(a\attach b) \circ f_a$.
\end{df}

Sheaf morphisms play an important role in models of signal processing because they transform global sections of one sheaf into those of another.

\begin{lem}
\label{lem:induced}
A sheaf morphism $f:\shf{F}\to\shf{G}$ induces a linear map $f_*:\Gamma \shf{F} \to \Gamma \shf{G}$.
\end{lem}

\section{Network monitoring}

The flow of contaminants through a water distribution network is usually assessed by measurements taken at specific locations.  These measurements are fundamentally local procedures, since they only collect a sample in the immediate vicinity of a location and time.  Because of this, we construct a sheaf that models the concentration of water-borne contaminants in a network of channels whose connection graph is represented by a directed graph $X$.  The edges of $X$ are labelled with a real-valued function $R$ which represents the volume flow rate of water along that edge.  In order to be physically reasonable, $R$ must satisfy a conservation law: at each vertex, the total flow rate in must equal the total flow rate out.  Additionally, we assume that any contaminant present in the inputs to a vertex is perfectly mixed at the vertex and flows out with the same concentration along all outgoing edges.

\begin{df}
The \emph{concentration sheaf} $\shf{C}$ on $X$ with flow rate $R$ is given by
\begin{enumerate}
\item $\shf{C}(v)=\mathbb{R}^n$ over each vertex $v$ with in-degree $n$, and
\item $\shf{C}(e)=\mathbb{R}$ over each edge $e$.
\end{enumerate}
The restriction from a vertex $v$ to the $i$-th incoming edge $e_i^+$ is given by the projection $\pr_i$ onto the $i$-th component of $\shf{C}(v)$.  Because of the perfect mixing assumptions, the restriction from a vertex $v$ to any outgoing edge $e^-$ is given by the formula
\begin{equation}
\label{eq:conc_rest}
\left(\shf{C}(v\attach e^-)\right)(c_1,\dotsc,c_n)=\frac{\sum_{j=1}^n c_j R(e^+_j)}{\sum_{j=1}^n R(e^+_j)}.
\end{equation}
\end{df}

Global sections of a concentration sheaf represent self-consistent values of contaminant concentration over the entire network.

\begin{eg}
Consider a network that distributes water from a main to several delivery points.  The concentration sheaf of this network is given by the diagram
\begin{equation*}
\xymatrix{
&&&\mathbb{R}&\mathbb{R}\ar[l]^{\id}\\
\mathbb{R}\ar[r]^{\id}&\mathbb{R}&\mathbb{R}\ar[l]_{\id}\ar[ur]^{\id}\ar[dr]_{\id}&\vdots&\\
&&&\mathbb{R}&\mathbb{R}\ar[l]^{\id}\\
}
\end{equation*}
where water flows from left-to-right.

By inspection, the space of sections is determined uniquely by measurement at any one of the delivery points.  Therefore, the contaminant concentration at the source can be recovered from a measurement at any of the delivery points.  Because of this reason, it is easy to ensure safe drinking water in a sealed distribution system (and track contaminations to their source) using endpoint checks.
\end{eg}

If concentration measurements are taken at vertices in the network $X$, they ought to be self-consistent.  This suggests that the process of taking measurements in a network is represented by a sheaf morphism.  

\begin{df}
A \emph{sampling morphism} of a concentration sheaf $\shf{C}$ on $X$ is a morphism $m:\shf{C}\to \shf{M}$ to a sheaf on $X$ which is the zero map on each edge and is surjective on the stalks over edges.  The assignment $\shf{A}(a)=\ker m_a$ defines another sheaf, called the \emph{ambiguity sheaf} over $X$. The restrictions of the ambiguity sheaf are given by the restrictions in $\shf{C}$, but restricted to the stalks of $\shf{A}$.
\end{df}

The consistent specification of concentrations throughout the network is a global section of the concentration sheaf $\shf{C}$, and a collection of measurements is a global section of the sheaf $\shf{M}$.  If the map induced on global sections by a sampling morphism (Lemma \ref{lem:induced}) is one-to-one, each global section of $\shf{M}$ corresponds to at most one global section of $\shf{C}$.  In this case, the samples contain sufficient information to reconstruct the concentrations over the entire network.

Although the induced map on global sections can be examined directly, it is helpful to have an alternate way to ensure that sufficiently many samples have been collected.  This characterization is provided by the following sampling theorem, which is valid for all networks (see \cite{robinson_2013} for more examples). 

\begin{thm}
If the ambiguity sheaf for a sampling morphism $m:\shf{C}\to\shf{M}$ has no nontrivial global sections, then the set of measurements completely specify all concentrations in the network.
\end{thm}

\begin{IEEEproof}(Sketch)
A well-known algebraic construction called the Snake Lemma shows that the kernel of the induced map $m_*:\Gamma\shf{C} \to \Gamma\shf{M}$ is given by the global sections of the ambiguity sheaf.  Therefore, if the ambiguity sheaf has no nontrivial global sections, then the induced map $m_*$ is one-to-one.
\end{IEEEproof}

\begin{eg}
Consider the concentration sheaf given by the diagram
\begin{equation*}
\xymatrix{
\mathbb{R}\ar[r]^{\id}&\mathbb{R}&&&\\
&\vdots& \mathbb{R}^n \ar[ul]^{\pr_1}\ar[dl]_{\pr_n}\ar[r]^M&\mathbb{R}&\mathbb{R}\ar[l]^{\id}\\
\mathbb{R}\ar[r]^{\id}&\mathbb{R}&&&\\
}
\end{equation*}
in which the water moves left-to-right, and $M$ is given by the $1\times n$ matrix representation of Equation \eqref{eq:conc_rest}.  The sheaf associated to measuring the concentrations at the delivery point has the diagram
\begin{equation*}
\xymatrix{
0\ar[r]^{\id}&0&&\\
&\vdots& 0 \ar[ul]^{\id}\ar[dl]_{\id}\ar[r]^{\id}&0&\mathbb{R}\ar[l]^{0}\\
0\ar[r]^{\id}&0&&&\\
}
\end{equation*}
The two sheaves are linked by a morphism, which is zero on all faces except the rightmost vertex, where it is an identity.  On the other hand, the ambiguity sheaf has the diagram
\begin{equation*}
\xymatrix{
\mathbb{R}\ar[r]^{\id}&\mathbb{R}&&\\
&\vdots& \mathbb{R}^n \ar[ul]^{\pr_1}\ar[dl]_{\pr_n}\ar[r]^{M}&\mathbb{R}&0\ar[l]^{0}\\
\mathbb{R}\ar[r]^{\id}&\mathbb{R}&&&\\
}
\end{equation*}
The space of global sections for the ambiguity sheaf (also the kernel of the induced map of the sampling morphism) has dimension $n-1$, since global sections are parameterized by the subspace of $\mathbb{R}^n$ on which $M$ vanishes.  This means that it is impossible to determine which input branch is responsible for the presence of a contaminant on the output.
\end{eg}

This example indicates why environmental management is a difficult problem.  Endpoint measurement cannot assign blame to polluters in water collection networks, but can be so used in water distribution networks.   Since industrial areas are often along rivers with the topology of a water collection network, their nominal cooperation ensures compliance with environmental regulations.

\section{Linear translation-invariant filters}
\begin{figure}
\begin{center}
\includegraphics[height=1in]{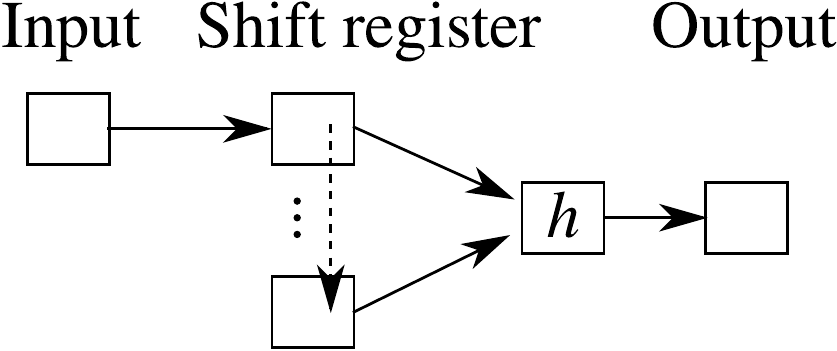}
\caption{Schematic diagram of a FIR LTI filter}
\label{fig:filter_wiring}
\end{center}
\end{figure}

Discrete linear translation-invariant filters are prevalent in modern signal processing systems, because they are particularly easy to implement.  If the impulse response of a filter is $N$ timesteps long, then the filter can be constructed using an $N$-word shift register as shown in Figure \ref{fig:filter_wiring}.  

\begin{df}
Let $l^\infty(V)$ be the vector space of bounded sequences of elements of a vector space $V$.  A \emph{discrete linear translation-invariant filter} (LTI filter) is a linear operator $F:l^\infty(V) \to l^\infty(V)$ given by
\begin{equation*}
\left(F(x)\right)_n = \sum_{k=-\infty}^\infty h(k) x_{n-k},
\end{equation*}
where $h:\mathbb{Z}\to \mathbb{R}$ is called the \emph{impulse response} of the filter.  If the support of $h$ is a finite set, we say that $F$ has \emph{finite impulse response} (FIR).
\end{df}

\begin{thm}
\label{thm:fir_encoding}
Every FIR LTI filter $F$ arises as the composition of linear maps $F=\lambda_*\circ p_*^{-1}:\Gamma \shf{S}_1 \to \Gamma \shf{S}_3$ induced on global sections by a pair of sheaf morphisms
\begin{equation*}
\xymatrix{
\shf{S}_1 & \shf{S}_2 \ar[r]^\lambda \ar[l]_p & \shf{S}_3.\\
}
\end{equation*}
In this diagram, the invertible linear map $p_*:\Gamma \shf{S}_2\to\Gamma \shf{S}_1$ is induced by $p$, and the map induced by $\lambda$ is $\lambda_*:\Gamma \shf{S}_2\to\Gamma\shf{S}_3$.
\end{thm}

This theorem has a clear interpretation in terms of the typical hardware implementation shown in Figure \ref{fig:filter_wiring}.  The global sections of $\shf{S}_1$ are precisely the possible input sequences, the global sections of $\shf{S}_3$ correspond to the output sequences, and the global sections of $\shf{S}_2$ correspond to the contents of the shift register.  The proof of the theorem is a rather explicit construction, which outlines the evolution of these three timeseries.

\begin{IEEEproof}
Let $\shf{S}_1$ and $\shf{S}_3$ both be copies of $\shf{V}^{(1)}$, the 1-term grouping sheaf.  The global sections of this sheaf are simply timeseries of data present in the input and output registers of the filter.  Suppose that the impulse response of $F$ is of length $N$ and let $\shf{S}_2=\shf{V}^{(N)}$ be the sheaf that represents the contents of the filter's shift register.

We construct the morphisms $p$ and $\lambda$ as the vertical arrows in the following diagram
\begin{equation*}
\xymatrix{
\ar[r]&0&V\ar[r]\ar[l]&0&\ar[l]\\
\ar[r]&V^{N-1}\ar[u]\ar[d]&V^N\ar[u]^{\pr_N}\ar[d]^L\ar[r]^{\sigma_+}\ar[l]^{\sigma_-}&V^{N-1}\ar[u]\ar[d]&\ar[l]\\
\ar[r]&0&V\ar[r]\ar[l]&0&\ar[l]\\
}
\end{equation*}
in which the first row is the diagram of $\shf{S}_1$, the second row represents $\shf{S}_2$, and the third row represents $\shf{S}_3$.  The maps $\pr_N$ and $L$ are given by the formulas
\begin{equation*}
\pr_N(x_1,\dotsc,x_N)=x_N,\text{ and }L(x_1,\dotsc,x_N)=\sum_{k=0}^{N-1} h(k) x_{N-k}.
\end{equation*}

It remains to verify that the induced map $\lambda_* \circ p_*^{-1}$ is the same as $F$.  Consider a sequence $x=(\dotsc,x_0,x_1,\dotsc)\in l^\infty(V)$, which can be encoded as a global section $s_1$ of $\shf{S}_1$ by specifying that the value of the section over the vertex $n$ is $x_n$, and the value over every edge is $0$.  This corresponds to a global section $s_2$ of $\shf{S}_2$, in which the value over each vertex is a sequence of $N$ consecutive terms of $x$, and the value over each edge is a sequence of $N-1$ consecutive terms.  Observe that the map $p_*$ induced by the morphism $p$ satisfies the equation $p_*(s_2)=s_1$, which moreover is one-to-one.  Therefore, with a slight abuse of notation, we write $p_*^{-1}(x) = s_2$.  The map $\lambda_*$ induced by the morphism $\lambda$ clearly computes weighted sums of adjacent blocks of $N$ consecutive terms of $x$, whence
\begin{eqnarray*}
L(x_1,\dotsc,x_N)&=&\sum_{k=0}^{N-1} h(k) x_{N-k}\\
&=&\sum_{k=-\infty}^\infty h(k) x_{N-k}\\
&=&\left(F(x)\right)_N.\\
\end{eqnarray*}
\end{IEEEproof}

\begin{eg}
Consider the FIR LTI filter whose impulse response is zero except for three consecutive terms, all equal to $1/3$.  If this filter is presented with the input sequence $\dotsc,1,1,9,2,\dotsc$, it will produce the output sequence $\dotsc,2.7,2.3,3.7,4,\dotsc$.  The encoding described by Theorem \ref{thm:fir_encoding} can be organized into the diagram
\begin{equation*}
\xymatrix{
0&2\ar[r]\ar[l]&0&9\ar[r]\ar[l]&0\\
(9,2)\ar[u]\ar[d]&(1,9,2)\ar[u]\ar[d]\ar[r]\ar[l]&(1,9)\ar[u]\ar[d]&(1,1,9)\ar[u]\ar[d]\ar[r]\ar[l]&(1,1)\ar[u]\ar[d]\\
0&4\ar[r]\ar[l]&0&3.7\ar[r]\ar[l]&0\\
}
\end{equation*}
\end{eg}

\begin{rem}
The benefit of using sheaf morphisms to describe filters is that they can treat a number of additional cases.  For instance, each of the following are straightforward generalizations, requiring no additional theoretical work to construct:
\begin{enumerate}
\item Infinite impulse response filters can be constructed simply by extending the definition of $\shf{V}^{(N)}$ to treat spaces of sequences instead of finite-dimensional vectors.
\item Nonlinear, block processing filters can be constructed by modifying the component maps of the morphism $\lambda$ to be nonlinear functions.  For instance, constant false alarm rate (CFAR) detectors can be encoded in this way.
\item If $G$ is a finitely-generated group that acts on $X$, then $G$-equivariant simplicial maps can be used to generalize $\shf{V}^{(N)}$ to other simplicial complexes $X$.  This permits extensions of Theorem \ref{thm:fir_encoding} to treat images, video, and more complex discrete datasets.
\end{enumerate}
\end{rem}

\section{Conclusion and discussion}
Sheaf theory offers a more general framework for extending key ideas (filtering and sampling) used in signal processing.  Sheaf theoretic examples of engineering problems often map directly to existing implementations.  Because of this, the use of sheaf theoretic tools requires minimal change in engineering perspective.

In the near term, we can expect that sheaf morphisms will come to the forefront of their respective engineering applications.  A more careful analysis of the morphisms that describe network measurements will probably result in tighter bounds on sampling requirements in general spaces with more general data than is possible with traditional methods.  Similarly, it is easy to extend the construction in Theorem \ref{thm:fir_encoding} to sheaves over different spaces.  However, the capabilities of the resulting \emph{topological filters} are essentially unexplored.  

\section*{ACKNOWLEDGEMENT}
This work was partly supported under Federal Contract No. FA9550-09-1-0643.

\bibliographystyle{IEEEtran}
\bibliography{globalsip2013_bib}

\begin{thebibliography}{1}
\providecommand{\url}[1]{#1}
\csname url@samestyle\endcsname
\providecommand{\newblock}{\relax}
\providecommand{\bibinfo}[2]{#2}
\providecommand{\BIBentrySTDinterwordspacing}{\spaceskip=0pt\relax}
\providecommand{\BIBentryALTinterwordstretchfactor}{4}
\providecommand{\BIBentryALTinterwordspacing}{\spaceskip=\fontdimen2\font plus
\BIBentryALTinterwordstretchfactor\fontdimen3\font minus
  \fontdimen4\font\relax}
\providecommand{\BIBforeignlanguage}[2]{{%
\expandafter\ifx\csname l@#1\endcsname\relax
\typeout{** WARNING: IEEEtran.bst: No hyphenation pattern has been}%
\typeout{** loaded for the language `#1'. Using the pattern for}%
\typeout{** the default language instead.}%
\else
\language=\csname l@#1\endcsname
\fi
#2}}
\providecommand{\BIBdecl}{\relax}
\BIBdecl

\bibitem{Kashiwara_1990}
M.~Kashiwara and P.~Schapira, \emph{Sheaves on manifolds}.\hskip 1em plus 0.5em
  minus 0.4em\relax Springer, 1990.

\bibitem{Baclawski_1975}
K.~Bac{\l}awski, ``Whitney numbers of geometric lattices,'' \emph{Adv. in
  Math.}, vol.~16, pp. 125--138, 1975.

\bibitem{Shepard_1980}
A.~Shepard, ``A cellular description of the derived category of a stratified
  space,'' Ph.D. dissertation, Brown University, 1980.

\bibitem{Lilius_1993}
J.~Lilius, ``Sheaf semantics for {Petri} nets,'' Helsinki University of
  Technology, Digital Systems Laboratory, Tech. Rep., 1993.

\bibitem{ghrist_2011}
R.~Ghrist and Y.~Hiraoka, ``Applications of sheaf cohomology and exact
  sequences to network coding,'' \emph{preprint}, 2011.

\bibitem{RobinsonLogic}
M.~Robinson, ``Asynchronous logic circuits and sheaf obstructions,''
  \emph{Electronic Notes in Theoretical Computer Science}, pp. 159--177, 2012.

\bibitem{robinson_2013}
------, ``The {Nyquist} theorem for cellular sheaves,'' in \emph{Sampling
  Theory and Applications (SampTA)}, 2013.

\bibitem{Curry}
J.~Curry, ``Sheaves, cosheaves and applications, {\tt arxiv:1303.3255},'' 2013.

\end{thebibliography}

\end{document}